\documentclass{article}
\usepackage{amsfonts}

\newtheorem{theorem}{Theorem}

\newtheorem{definition}[theorem]{Definition}
\newtheorem{example}[theorem]{Example}
\newtheorem{remark}[theorem]{Remark}

\def\QED{\quad\blackslug\lower 8.5pt\null}

\begin{document}

\begin{center}
{\large \bf Local Equivalence of Sacksteder and Bourgain Hypersurfaces}
\end{center}

\vspace*{3mm}

\begin{center}
{\large
 Maks A. Akivis and  Vladislav V. Goldberg}
\end{center}

\vspace*{3mm}

{\footnotesize
\noindent
{\em $2000$ Mathematics Subject Classification}.  
Primary 53A20, Secondary 14M99.

\textbf{Abstract.}
Finding examples of tangentially degenerate submanifolds
(submanifolds with degenerate Gauss mappings)
in an Euclidean space $R^4$ that are  noncylindrical
and  without singularities is an important problem of differential
geometry. The first example of such a hypersurface was constructed
by Sacksteder in 1960. In 1995
Wu published an example of a noncylindrical tangentially degenerate algebraic
hypersurface in  $R^4$  whose Gauss mapping is of rank 2 and
which  is also without singularities.
This example was constructed (but not published) by Bourgain.
 \par  In this paper, the authors analyze Bourgain's example, prove that,
  as was the case for the Sacksteder hypersurface,
singular points of the Bourgain hypersurface are located in the hyperplane
at infinity  of the space $R^4$, and  these
   two hypersurfaces are locally equivalent.
}

\setcounter{equation}{0}

\vspace*{5mm}

\textbf{1.}
It is important to find examples of tangentially degenerate submanifolds
in order to understand the theory of such manifolds.
These examples prove the existence of tangentially degenerate
submanifolds and help to illustrate the theory.
   The first known example of a tangentially degenerate hypersurface
of rank 2 without singularities in $R^4$ was constructed by Sacksteder [S 60].
This example was examined from the differential geometry point of view
by Akivis in [A 87]. In particular, Akivis proved that the Sacksteder
 hypersurface has no singularities since they ``went to infinity''.
  In the same paper, Akivis presented
a series of examples generalizing Sacksteder's example in $R^4$,
constructed a new series of examples of three-dimensional hypersurfaces
$V^3 \subset P^n (\mathbb{R})$ of rank 2 whose focal surfaces are imaginary, and
proved existence of hypersurfaces of this kind.
Note that more examples of
tangentially degenerate submanifolds without singularities
can be found in [I 98, 99a, 99b]. The examples are essentially
based on classical Cartan's hypersurfaces (see [Ca 39]).

Mori [M 94] claims that he constructed ``a one-parameter family
of complete nonruled deformable hypersurfaces in $R^4$
with rank $r = 2$ almost everywhere''. However, it follows
immediately from his formulas that the hypersurfaces of his family
are not only ruled hypersurfaces---they are cylinders.

Also much progress on the study of tangentially degenerate
submanifolds over the complex numbers has been made in
[GH 79], [L 99], and [AGL 00]. In these papers and in the papers
[FW 95], [W 95], and [WZ 99] one can find more examples
of tangentially degenerate submanifolds over the complex numbers.

Recently Wu [W 95] published an example of a noncylindrical
tangentially degenerate algebraic hypersurface in an Euclidean
space $R^4$ which has a degenerate Gauss mapping but does not
have singularities. This example was constructed (but not
published) by Bourgain (see also [I 98, 99a, 99b]).
In the present paper, we
investigate  Bourgain's example from the point of view of
the paper [A 87] (see also Section 4.7 of our book [AG 93]).
In particular, we prove that,  as was the case for the Sacksteder  hypersurface,
the Bourgain  hypersurface  has no singularities since they ``went to
infinity''. Namely this analysis suggested an idea that Bourgain's and Sacksteder's
examples must be equivalent. Moreover, this analysis
showed that a hypersurface constructed in these examples is
torsal, i.e., it is stratified into a one-parameter family of
plane pencils of straight lines.

 In addition, at the end of our paper we prove that  the examples of  Bourgain
 and  Sacksteder are locally equivalent.

\textbf{2.} In Cartesian coordinates $x_1, x_2, x_3, x_4$ of
the Euclidean space $R^4$, the equation of
the Bourgain hypersurface $B$ is
\begin{equation}\label{1}
  x_1 x_4^2 + x_2 (x_4 - 1) + x_3 (x_4 - 2) = 0
\end{equation}
(see [W 95] or [I 98, 99a, 99b]). Equation (1) can be written
in the form
\begin{equation}\label{2}
  x_1 x_4^2 + (x_2 + x_3) x_4 - (x_2 + 2 x_3) = 0.
\end{equation}
Make in (2) the following admissible  change of Cartesian coordinates:
$$
 x_2 + x_3 \rightarrow x_2, \;\;
 x_2 + 2 x_3 \rightarrow x_3.
 $$
 Then equation (2) becomes
\begin{equation}\label{3}
  x_1 x_4^2 + x_2 x_4 - x_3 = 0.
\end{equation}

Introduce homogeneous coordinates in $R^4$
by setting $ x_i = \displaystyle\frac{z_i}{z_0}, \;
i = 1, 2, 3, 4$. Then equation (3) takes the form
\begin{equation}\label{4}
 f = z_1 z_4^2 + z_0 z_2 z_4 - z_0^2  z_3 = 0.
\end{equation}
Equation (4) defines a cubic hypersurface $F$ in the space
$\overline{R}^4 = R^4 \cup P^3_\infty$
which is an enlarged space $R^4$, i.e., it is
the  space $R^4$ enlarged by the hyperplane at
infinity $ P^3_\infty$ (whose equation is $z_0 = 0$).

Denote by $A_\alpha, \; \alpha = 0, 1, 2, 3, 4,$
fixed basis points of the space $\overline{R}^4$.
Suppose that these points have constant normalizations,
i.e., that $dA_\alpha = 0$. An arbitrary point $z
\in  \overline{R}^4$ can be written in the form
$z = \sum_\alpha z_\alpha A_\alpha$.
We will take a proper point of the space
$\overline{R}^4$ as the point  $A_0$,
and take points at infinity as the
 points  $A_1, A_2, A_3, A_4$.

 Equation (4) shows that the proper straight
 line $A_0 \wedge A_4$ defined by the
 equations $z_1 = z_2 = z_3 = 0$ and the plane at
 infinity defined by the equations
 $z_0 = z_4 = 0$ belong
 to the hypersurface $F$ defined by
 equation (4).

 We write the equations of the hypersurface $F$
 in a parametric form. To this end, we set
 $$
 z_0 = 1, \;\;z_4 = p, \;\; z_1 = u,\;\; z_3 = p v.
 $$
 Then it follows from (4) that
 $$
 z_2 = v - p u.
 $$
 This implies that an arbitrary point $z \in F$
 can be written as
\begin{equation}\label{5}
  z = A_0 + u A_1 + v A_2 + p(A_4 - u A_2 + v A_3).
\end{equation}
The parameters $p,\; u, \;v$ are independent nonhomogeneous
parameters on  the hypersurface $F$.

\textbf{3.} Let us find singular points of  the hypersurface $F$.
Such points are defined by the equations
$\displaystyle\frac{\partial
f}{\partial z_\alpha} = 0$. It follows from (4) that
\begin{equation}\label{6}
\renewcommand{\arraystretch}{2.1}
\left\{
\begin{array}{ll}
\displaystyle \frac{\partial f}{\partial z_0}
= z_2 z_4 - 2 z_0 z_3,\\
 \displaystyle \frac{\partial f}{\partial z_1} = z_4^2,\;\;\;
\displaystyle \frac{\partial f}{\partial z_2} = z_0 z_4, \;\;\;
\displaystyle \frac{\partial f}{\partial z_3} = - z_0^2,\\
\displaystyle \frac{\partial f}{\partial z_4} = 2 z_1 z_4 + z_0 z_2.
\end{array}
\right.
\renewcommand{\arraystretch}{1}
\end{equation}
All these derivatives vanish simultaneously if and only if $z_0 =
z_4 = 0$. Thus the 2-plane at infinity $\sigma = A_1 \wedge A_2
\wedge A_3$ is the locus of singular points of  the hypersurface $F$.

Consider a point $B_0 = A_0 + p A_4$ on the straight line
$ A_0 \wedge A_4$. By (4), to the point $B_0$ there corresponds
the straight line $a (p)$ in the 2-plane  at infinity
$\sigma$, and the equation of this straight line is
\begin{equation}\label{7}
  p^2 z_1 + p z_2 - z_3 = 0.
\end{equation}
The family of straight lines $a (p)$ depends of the parameter $p$,
and its envelope is the conic $C$ defined by the equation
\begin{equation}\label{8}
 z_2^2 + 4 z_1 z_3 = 0.
\end{equation}
The straight line $a (p)$ is tangent to  the conic $C$
at the point
\begin{equation}\label{9}
B_1 (p) = A_1 - 2 p A_2 - p^2 A_3.
\end{equation}
Equation (9) is a parametric equation of  the conic $C$. The point
\begin{equation}\label{10}
 \displaystyle \frac{d B_1}{d p} = - 2 (A_2 + p A_3)
\end{equation}
 belongs to the tangent line to  the conic $C$ at the point $B_1
 (p)$.

 Consider the 2-planes $\tau = B_0 \wedge B_1 \wedge  \displaystyle \frac{d B_1}{d
 p}$. Such 2-planes are completely determined by the location of
 the point $B_0$ on  the straight line
$ A_0 \wedge A_4$, and they form a one-parameter family.
All these 2-planes belong to the hypersurface $F$. In fact,
represent an arbitrary point $z$ of the 2-plane $\tau$ in the form
\begin{equation}\label{11}
\renewcommand{\arraystretch}{1.3}
\begin{array}{ll}
z &= \alpha B_0 + \beta B_1 - \displaystyle
\frac{1}{2} \gamma \displaystyle \frac{d B_1}{d p}\\
&= \alpha A_0 + \beta A_1 + (-2 p \beta + \gamma) A_2 +
(- p^2 \beta + p \gamma) A_3 + p \alpha A_4.
\end{array}
\renewcommand{\arraystretch}{1}
\end{equation}
  The coordinates of the point $z$ are
\begin{equation}\label{12}
 z_0 = \alpha, \;\;z_1 = \beta, \;\;z_2 = \gamma - 2 p \beta, \;\;
 z_3 = p (\gamma - p \beta), \;\;z_4 = p\alpha.
 \end{equation}
Substituting these values of the coordinates
into equation (4), one can see that equation (4)
is identically satisfied. Thus
the hypersurface $F$ is foliated into a one-parameter family
of 2-planes $\tau (p)
= B_0 \wedge B_1 \wedge  \displaystyle \frac{d B_1}{d  p}$.

In a 2-plane $\tau (p)$ consider a pencil of straight lines with center $B_1$.
The  straight lines of this pencil are
defined by the point $B_1$ and the point $B_2 = A_2 + p A_3 + q
(A_0 + p A_4)$. The straight lines $B_1 \wedge B_2$ depend on two
parameters $p$ and $q$. These lines belong to the 2-plane $\tau (p)$,
and along with this 2-plane they belong to
the hypersurface $F$. Thus they form a foliation
on the hypersurface $F$.

We prove that this foliation is a Monge-Amp\`{e}re foliation.
In the space  $\overline{R}^4$, we introduce the moving frame
formed by the points
\begin{equation}\label{13}
\renewcommand{\arraystretch}{1.3}
\left\{
\begin{array}{ll}
B_0 = A_0 + p A_4,\\
B_1 = A_1 - 2 p A_2 - p^2 A_3, \\
B_2 = A_2 + p A_3 + q A_0 + p q A_4, \\
B_3 = A_3,\\
B_4 = A_4.
\end{array}
\right.
\renewcommand{\arraystretch}{1}
\end{equation}
It is easy to prove that these points are linearly independent,
and the points $A_\alpha$ can be expressed in terms of the points
$B_\alpha$ as follows
\begin{equation}\label{14}
\renewcommand{\arraystretch}{1.3}
\left\{
\begin{array}{ll}
A_0 = B_0 - p B_4,\\
A_1 = B_1 + 2 p B_2 - p^2 B_3 - 2 p q B_0, \\
A_2 = B_2 - p B_3 - q B_0, \\
A_3 = B_3,\\
A_4 = B_4.
\end{array}
\right.
\renewcommand{\arraystretch}{1}
\end{equation}

Consider a displacement of the straight lines $B_1 \wedge B_2$
along the hypersurface $F$. Suppose that $Z$ is an arbitrary
point of this straight line,
\begin{equation}\label{15}
Z = B_1 + \lambda B_2.
\end{equation}
Differentiating (15) and taking into account (14) and
$d A_\alpha = 0$, we find that
\begin{equation}\label{16}
d Z \equiv (2 q d p + \lambda d q) B_0 + \lambda d p (B_3 + q B_4)
 \pmod{B_1, B_2}.
\end{equation}
It follows from relation (16) that
\begin{enumerate}
\item A tangent hyperplane to the hypersurface $F$
is spanned by the points $B_1, B_2, B_0$ and
$B_3 + q B_4$. This hyperplane is fixed when the point $Z$
moves along the straight line $B_1 \wedge B_2$.
Thus  the hypersurface $F$ is tangentially degenerate
of rank 2, and the straight lines $B_1 \wedge B_2$
form a Monge-Amp\`{e}re foliation  on $F$.

\item The system of equations
\begin{equation}\label{17}
\renewcommand{\arraystretch}{1.3}
\left\{
\begin{array}{ll}
2 q d p + \lambda d q & = 0, \\
\lambda d p           & = 0
\end{array}
\right.
\renewcommand{\arraystretch}{1}
\end{equation}
defines singular points on the straight line $B_1 \wedge B_2$,
and on the hypersurface $F$ it defines torses. The system of
equations (17) has a nontrivial solution with respect to $d p$ and
$d q$ if and only if its determinant vanishes: $\lambda^2 = 0$.
Hence by (15), a singular point on the straight line $B_1 \wedge B_2$
coincides with the point $B_1$. For $\lambda = 0$, system (17)
implies that $d p = 0$, i.e., $p = \mbox{{\rm const}}$.
Thus it follows from (9) that the point $B_1 \in C$ is fixed, and
as a result, the torse corresponding to this constant parameter
$p$ is a pencil of straight lines with the center $B_1$ located in the
2-plane $\tau (p) = B_0 \wedge B_1 \wedge B_2$.

\item All singular points of the hypersurface $F$ belong to
the conic $C \subset P^\infty$ defined by equation (8).
Thus if we consider  the hypersurface $F$  in an Euclidean
space $R^4$, then on $F$ there are no singular points
in a proper part of this space.

\item  The hypersurface $F$ considered in the proper part of
an Euclidean space is not a cylinder since its rectilinear
generators do not belong to a bundle of parallel straight
lines. A two-parameter family of  rectilinear
generators of $F$ decomposes into a one-parameter family
of plane pencils of parallel lines.
\end{enumerate}

\textbf{4.}
    No one of properties 1--4  characterizes Bourgain's
    hypersurfaces completely: they are necessary but not
    sufficient for these hypersurfaces. The following
  theorem   gives a   necessary and sufficient condition
    for a hypersurface to be of Bourgain's  type.

\begin{theorem}
Let $l$ be a proper straight line of an Euclidean space $R^4$
enlarged by the plane at infinity $P^3_\infty$, and let $C$
be a conic in the $2$-plane $\sigma$. Suppose that the straight line $l$ and
the conic $C$ are in a projective correspondence. Let $B_0 (p)$
and $B_1 (p)$ be two corresponding points of $l$ and $C$, and
let $\tau$ be the $2$-plane passing through the point $B_0$
and tangent to the conic $C$ at the point $B_1$. Then
\begin{itemize}
\item[{\rm (a)}] when
the point $B_0$ is moving along the straight line $l$, the plane
$\tau$ describes a Bourgain hypersurface, and
\item[{\rm (b)}] any  Bourgain hypersurface satisfies the above construction.
\end{itemize}
\end{theorem}

{\sf Proof}.
The necessity (b) of the theorem hypotheses follows from our
previous considerations.
We prove the sufficiency (a)  of these hypotheses.
Take a fixed frame $\{A_u\}, \, u = 0, 1, 2, 3, 4,$
in the  space $R^4$
enlarged by the plane at infinity $P^3_\infty$ as follows:
 its point $A_0$ belongs to $l$, the point
$A_4$ is the point at infinity of $l$, and  the points
$A_1, A_2$, and $A_3$ are located at the $2$-plane at infinity
$\sigma$ in such a way that a parametric equation of
the straight line $l$ is  $B_0 = A_0 + p A_4$, and
the equation of $C$ has the form (9). The plane $\tau$
is defined  by  the points
$B_0, B_1$, and $\frac{d B_1}{d p}$. The parametric
equations of this plane have the  form (12). Excluding
the parameters $\alpha,\beta, \gamma$, and $p$  from
these equations, we will return to the cubic equation (4)
defining the Bourgain hypersurface $B$ in homogeneous coordinates.
\rule{3mm}{3mm}

The method of construction of the Bourgain  hypersurface
used in the proof of Theorem 1 goes back to
the classical methods of projective geometry developed by Steiner
[St 32] and Reye [R 66].

\textbf{5.} In conclusion we prove the following theorem.

\begin{theorem}
The Sacksteder hypersurface $S$ and the Bourgain hypersurface $B$
are locally equivalent, and the former is the
standard covering of the latter.
\end{theorem}

{\sf Proof.}
In an Euclidean space $R^4$, in
 Cartesian coordinates $x_1, x_2, x_3, x_4$,
 the equation of the  Sacksteder hypersurface $S$ (see [S 60])
 has the form
\begin{equation}\label{18}
  x_4 = x_1 \cos x_3 + x_2 \sin x_3.
\end{equation}
The right-hand side of this equation is a function on
the manifold $M^3 = \mathbb{R}^2 \times S^1$ since
the variable $x_3$ is cyclic. Equation (18) defines
a hypersurface on the manifold $M^3 \times S^1$.
The circumference $S^1 = \mathbb{R}/2 \pi \mathbb{Z}$ has
a natural projective structure of $P^1$. In the homogeneous
coordinates $x_3 = \frac{u}{v}$, the mapping $S^1 \rightarrow
P^1$, can be written as $x^3 \rightarrow (u, v)$.
By removing the  point $\{v = 0\}$ from $S^1$, we obtain
a 1-to-1 correspondence
\begin{equation}\label{19}
 S^1 - \{v = 0\} \longleftrightarrow \mathbb{R}^1.
\end{equation}
Now we can consider the  Sacksteder hypersurface $S$
in $R^4$ or, if we enlarge $R^4$ by the plane at infinity
$P^3_{\infty}$, in the space $P^4$.

Next we show how by applying the mapping $S^1 \rightarrow P^1$,
we can transform equation (18) of the  Sacksteder hypersurface $S$
into equation(4) of the  Bourgain hypersurface $B$.
We write this mapping in the form
\begin{equation}\label{20}
x_3 = 2 \arctan \frac{u}{v}, \;\; \frac{u}{v} \in R, \; |x_3| <
\pi.
\end{equation}
It follows from (20) that
\begin{equation}\label{21}
\renewcommand{\arraystretch}{1.3}
\left\{
\begin{array}{ll}
\displaystyle\frac{u}{v} = \tan \displaystyle\frac{x_3}{2}, \\
\cos x_3 = \displaystyle\frac{1 - \tan^2  \displaystyle\frac{x_3}{2}}{1 + \tan^2
 \displaystyle\frac{x_3}{2}} = \displaystyle\frac{v^2 - u^2}{v^2 + u^2}, \\
\sin x_3 = \displaystyle\frac{2 \tan  \displaystyle\frac{x_3}{2}}{1 + \tan^2
 \displaystyle\frac{x_3}{2}} = \displaystyle\frac{2 u v}{v^2 + u^2}.
\end{array}
\right.
\renewcommand{\arraystretch}{1}
\end{equation}
Substituting these expressions into equation (18), we find that
$$
x_4 (u^2 + v^2) = x_1 (v^2 - u^2) + 2 x_2 uv,
$$
i.e.,
\begin{equation}\label{22}
  (x_4 + x_1) u^2 + (x_4 - x_1) v^2 - 2 x_2 uv = 0.
\end{equation}
Make a change of variables
$$
  z_1 = x_4  - x_1, \;\; z_2 = -2 x_2, \;\; z_3 = x_1 + x_4, \;\;
  z_0 = u, \;\; z_4 = v.
$$
As a result, we reduce equation (22) to equation (4).
It follows that the Sacksteder hypersurface $S$ defined by equation (18)
is locally equivalent to the Bourgain hypersurface
defined by equation (4).

Note also that if the cyclic parameter $x_3$ changes on the entire
real axis $R$, then we obtain the standard
covering of the Bourgain  hypersurface $B$ by means of
the Sacksteder hypersurface $S$.
\rule{3mm}{3mm}

\noindent {\em Authors' addresses}:\\

\noindent
\begin{tabular}{ll}
M.~A. Akivis &V.~V. Goldberg\\
Department of Mathematics &Department of Mathematical Sciences\\
Jerusalem College of Technology---Mahon Lev &  New
Jersey Institute of Technology \\
Havaad Haleumi St., P. O. B. 16031 & University Heights \\
 Jerusalem 91160, Israel &  Newark, N.J. 07102, U.S.A. \\
 & \\
 E-mail address: akivis@avoda.jct.ac.il & E-mail address:
 vlgold@m.njit.edu
 \end{tabular}
\end{document}